\documentclass[a4paper,11pt]{amsart}
\usepackage{amsmath,amsthm}
\newtheorem{theorem}{Theorem}[section]
\newtheorem{cor}[theorem]{Corollary}
\newtheorem{lemma}[theorem]{Lemma}
\newtheorem{prop}[theorem]{Proposition}
\newtheorem{scholium}[theorem]{Scholium}
\theoremstyle{remark}
\newtheorem{remark}[theorem]{Remark}
\newtheorem{example}[theorem]{Example}
\theoremstyle{definition}
\newtheorem{definition}[theorem]{Definition}
\numberwithin{equation}{section}
\DeclareMathOperator{\Isom}{Isom}
\DeclareMathOperator{\cv}{cov}
\DeclareMathOperator{\clsp}{\overline{span}}
\newsymbol\rtimes 216F
\newcommand{\abs}[1]{\lvert#1\rvert}
\newcommand{\norm}[1]{\lVert#1\rVert}
\newcommand{\cov}{\text{$C^*_{{\cv}}$}}
\newcommand{\cross}{\rtimes}
\newcommand{\bpp}{B_P \textstyle{\cross_{\tau, E}} P}
\newcommand{\rankone}[2]{#1 \otimes\overline{#2}}
\newcommand{\field}[1]{\mathbb{#1}}
\newcommand{\NN}{\field{N}}
\newcommand{\RR}{\field{R}}
\newcommand{\TT}{\field{T}}
\newcommand{\Bb}{{\mathcal B}}
\newcommand{\Cc}{{\mathcal C}}
\newcommand{\Ee}{{\mathcal E}}
\newcommand{\Gg}{{\mathcal G}}
\newcommand{\Hh}{{\mathcal H}}
\newcommand{\Kk}{{\mathcal K}}
\newcommand{\Oo}{{\mathcal O}}
\newcommand{\Pp}{{\mathcal P}}
\newcommand{\Tt}{{\mathcal T}}
\begin{document}
%
%
\title[compactly-aligned discrete product systems]
{Compactly-aligned discrete product systems, \\
and generalizations of ${\mathcal O}_\infty$}
\author[Neal J. Fowler]{Neal J. Fowler}
\address{Department of Mathematics  \\
      University of Newcastle\\  NSW  2308\\ AUSTRALIA}
\email{neal@math.newcastle.edu.au}
\date{September 26, 1998}
\begin{abstract}
The universal $C^*$-algebras of discrete product systems
generalize the Toeplitz-Cuntz algebras
and the Toeplitz algebras of discrete semigroups.
We consider a semigroup $P$ which is quasi-lattice ordered
in the sense of Nica,
and, for a product system $p:E\to P$,
we study those representations of $E$,
called covariant, which respect the lattice structure of $P$.
We identify a class of product systems, which we call
compactly aligned, for which there is a purely $C^*$-algebraic
characterization of covariance,
and study the algebra $C^*_{\cv}(P,E)$ which is universal
for covariant representations of $E$.
Our main theorem is a characterization of the faithful representations
of $\cov(P,E)$ when $P$ is the positive cone of a free product of totally-ordered
amenable groups.
\end{abstract}
\maketitle
%
%
\section*{Introduction}

The study of $C^*$-algebras generated by isometries has a rich history, beginning
with Coburn's theorem that any two algebras $C^*(U)$ and $C^*(V)$ generated by nonunitary
isometries $U$ and $V$ are canonically isomorphic \cite{coburn}.
Cuntz generalized Coburn's result by studying the $C^*$-algebra generated by a
{\em Toeplitz-Cuntz} family $\{S_1, \dots, S_n\}$ of isometries with mutually orthogonal ranges;
this led to his well-known analysis of the Cuntz algebra $\Oo_n$
and its Toeplitz extension $\Tt\Oo_n$ \cite{cun,cun2}.

One can also view Coburn's algebra as the Toeplitz algebra of $\NN$, and it is natural
to seek versions of his theorem for other semigroups.
A modern approach is to take a semigroup $P$, construct a $C^*$-algebra $C^*(P)$
which is universal for representations $V:P\to\Isom(\Hh)$
(called {\em isometric representations\/}),
and ask when the integrated form $V_*:C^*(P)\to\Bb(\Hh)$ is faithful;
in this language Coburn's theorem says that $V_*:C^*(\NN)\to\Bb(\Hh)$
is faithful if and only if none of the isometries $V_n$ is unitary.
Douglas generalized this result by showing that one can replace $\NN$
with the positive cone $\Gamma^+$ of any countable subgroup $\Gamma$ of $\RR$
\cite{douglas}, and Murphy later extended this to handle any
totally-ordered abelian group $\Gamma$ \cite{murjot}.

Here we combine the approaches of Cuntz and Douglas/Murphy by studying the $C^*$-algebra
generated by an entire collection of Toeplitz-Cuntz families, one for each element
of a semigroup $P$.  The algebraic relations among the families are governed by
a {\em product system\/} $E$ over $P$, which, loosely speaking, is a collection
$\{E_s: s\in P\}$ of Hilbert spaces, together with an associative multiplication
which induces unitary isomorphisms $E_s \otimes E_t \to E_{st}$.
(Our $E$ is a discrete analogue of the continuous tensor product systems
which arose in Arveson's work on one-parameter semigroups of endomorphisms of $\Bb(\Hh)$
\cite{arv}.)   Every representation of $E$ gives a collection of Toeplitz-Cuntz families
indexed by $P$, and with an appropriate choice of $E$
we can realize any of the above Toeplitz algebras as the $C^*$-algebra
$C^*(P,E)$ which is universal for representations of $E$:
when $E$ is the trivial product system over $P$ we have $C^*(P,E) = C^*(P)$,
and by taking the $n$-dimensional product system over $\NN$ we have $C^*(P,E) = \Tt\Oo_n$.
(See \cite[\S1]{fowrae} for these and other examples.)

The $C^*$-algebra of a discrete product system was first studied by Dinh,
who established the simplicity of $C^*(\Gamma^+,E)$ for any product system $E$
over the positive cone $\Gamma^+$ of a countable dense subgroup $\Gamma$ of $\RR$,
so long as each fiber $E_s$ is infinite dimensional \cite{dinhjfa}.
In \cite{fowrae}, Fowler and Raeburn considered product systems over more general
semigroups, but focussed primarily on systems with finite-dimensional fibers.
Here we emphasize product systems which have infinite-dimensional fibers,
and one of our goals is to extend Dihn's theorem to Murphy's setting;
this is achieved in Scholium~\ref{scholium}.

Extending Murphy's result beyond the totally ordered case is problematic, in part because
for an arbitrary isometric representation $V$, the algebra $C^*(V)$ generated
by $\{V_s: s\in P\}$ need not be spanned by monomials of the form $V_rV_s^*$:
one cannot simplify a product $V_rV_s^*V_tV_u^*$ unless $s$ and $t$ are comparable.
In \cite{nica}, Nica made the important observation
that for a large class of partially ordered groups $(G,P)$, which he called
{\em quasi-lattice ordered} groups, this problem does not exist when $V$ is the
left regular representation of $P$ on $\ell^2(P)$: the range projections $V_sV_s^*$
respect the lattice structure of $P$ in the sense that
$V_sV_s^*V_tV_t^* = V_{s\vee t}V_{s\vee t}^*$,
and this additional algebraic relation allows the simplification
$V_rV_s^*V_tV_u^* = V_{rs^{-1}(s\vee t)}V_{ut^{-1}(s\vee t)}^*$.
Nica called isometric representations with this property {\em covariant\/},
and defined $C^*(G,P)$ as the $C^*$-algebra which is universal for
covariant isometric representations.

Since covariance is automatic when $(G,P)$ is a total order,
the class of algebras $C^*(G,P)$ includes those considered by Murphy;
in \cite{lacarae}, Laca and Raeburn used semigroup crossed product techniques to
characterize the faithful representations of $C^*(G,P)$, thus generalizing
Murphy's theorem.
Motivated by this work, Fowler and Raeburn defined a notion
of covariance for representations of a product system $E$ over a quasi-lattice order $P$,
defined a universal algebra $\cov(P,E)$, and characterized its faithful representations
when $E$ has finite-dimensional fibers \cite{fowrae}.

For systems with infinite-dimensional fibers the definition of covariance
involves infinite sums which converge $\sigma$-weakly but not in norm;
consequently, as observed in \cite[Remark~4.5]{fowrae},
one cannot be certain that every representation of $\cov(P,E)$ is the integrated
form of a covariant representation.
In Section~1 we give a concrete illustration of this pathology (Example~\ref{example}),
and then define a class of product systems,
called {\em compactly aligned\/}, for which there is no such problem:
there is a $C^*$-algebraic characterization of covariance
(Proposition~\ref{prop:covariance}),
and this assures that $\cov(P,E)$ behaves as a universal object should
(Proposition~\ref{prop:i_E covariant}).

Most of the known examples of discrete product systems are compactly aligned;
indeed, if $E$ has finite-dimensional fibers, or if $(G,P)$ is a total order,
then $E$ is compactly aligned, and the class of compactly-aligned product systems
is closed under free products (Proposition~\ref{prop:free products}).
One can also adapt Dihn's construction of von Neumann discrete product systems
to an arbitrary $(G,P)$, and these product systems are compactly aligned
as well (Proposition~\ref{prop:von Neumann}).

By showing that compactly-aligned product systems satisfy
the spanning condition required in \cite[Theorem~6.1]{fowrae},
we obtain a sufficient condition for faithfulness of a covariant representation
(Theorem~\ref{theorem:sufficient}), and in Corollary~\ref{cor:simplicity}
we use this condition to generalize a simplicity result of Laca and Raeburn.
To obtain a condition which is also necessary, we then specialize to product systems
over free products of totally-ordered amenable groups;
for such a product system $E$ our main result, Theorem~\ref{theorem:faithful},
characterizes the faithful representations of $\cov(P,E)$.
The proof is a modification of the one given in \cite[Theorem~5.1]{fowrae}:
the gauge action gives an expectation on $\cov(P,E)$,
and the key steps are to implement this expectation spatially,
and to prove faithfulness on the fixed-point algebra.  
The details are necessarily different; in particular, in performing the usual
norm estimates we need the free product structure to implement Cuntz's technique
of using an aperiodic sequence to kill the off-diagonal terms.

We would like to thank Iain Raeburn for many helpful discussions
during the preparation of this work.

\section{Compactly-aligned discrete product systems}

We begin by recalling the pertinent definitions from \cite{fowrae}.
Suppose $P$ is a countable discrete semigroup with identity $e$.
A {\em product system\/} over $P$ is a family $p:E\to P$
of nontrivial separable complex Hilbert spaces $E_t:= p^{-1}(t)$
which is endowed with an associative multiplication $E\times E\to E$
in such a way that $p$ is a semigroup homomorphism,
and such that for every $s,t\in P$ the map
$x\otimes y\in E_s\otimes E_t \mapsto xy\in E_{st}$ 
extends to a unitary isomorphism.
We also insist that $\dim E_e = 1$,
so that $E$ has an identity $\Omega$ \cite[Lemma~1.3]{fowrae}.

A {\em representation\/} of $E$ in a unital $C^*$-algebra $B$
is a map $\phi:E\to B$ which satisfies
\begin{itemize}
\item[(i)] $\phi(xy) = \phi(x)\phi(y)$ for every $x,y\in E$; and
\item[(ii)] $\phi(y)^*\phi(x) = \langle x,y \rangle 1$ if $p(x) = p(y)$.
\end{itemize}
It follows that the restriction of $\phi$ to each $E_t$
is linear and isometric, and that $\phi(\Omega) = 1$ \cite[Remarks~1.8]{fowrae}.
Condition (ii) implies that each $\phi(x)$ is a multiple of an isometry;
moreover, the image of an orthonormal basis for $E_t$ is a Toeplitz-Cuntz family,
and this collection of families indexed by $P$ generates $C^*(\phi(E))$.
When $B = \Bb(\Hh)$ for some separable Hilbert space $\Hh$, we call $\phi$ a
{\em representation of $E$ on $\Hh$\/}.

Given a representation $\phi:E\to B$, for each $s\in P$ there is a homomorphism
$\rho^\phi_s:\Kk(E_s)\to B$ which satisfies
\[
\rho^\phi_s(\rankone xy) = \phi(x)\phi(y)^*,\qquad x,y\in E_s;
\]
here $\rankone xy$ denotes the rank-one operator $z \mapsto \langle z,y \rangle x$.
It is easy to check that
\begin{equation}\label{eq:rho}
\rho^\phi_s(S)\phi(z) = \phi(Sz)
\qquad\text{for $S\in\Kk(E_s)$, $z\in E_s$,}
\end{equation}
and if $\sigma:B\to C$ is a homomorphism, then
\begin{equation}\label{eq:compose}
\rho^{\sigma\circ\phi}_s(S) = \sigma(\rho^\phi_s(S))
\qquad\text{for $S\in\Kk(E_s)$.}
\end{equation}

When $B = \Bb(\Hh)$, $\rho^\phi_s$ has a unique extension to $\Bb(E_s)$;
we write $\rho^\phi_s$ for this extension as well, and remark that
it is normal ($\sigma$-weakly continuous).
By \cite[Proposition~1.11]{fowrae}, there is also an associated semigroup homomorphism
$t \mapsto \alpha^\phi_t$
of $P$ into the ${}^*$-endomorphisms of $\Bb(\Hh)$, given by
\begin{equation}\label{eq:alpha}
\alpha^\phi_t(A) = \sum_{u\in\Bb} \phi(u)A\phi(u)^*
\qquad\text{for $A\in\Bb(\Hh)$,}
\end{equation}
where the sum converges strongly over any orthonormal basis $\Bb$ for $E_t$.

The following Lemma summarizes the properties of $\alpha^\phi$ and $\rho^\phi$ we shall require.

\begin{lemma}\label{lemma:alpharho}
If $s,t\in P$, $z\in E_s$, $S\in\Bb(E_s)$ and $A\in\Bb(\Hh)$, then

\textup{(1)} $\alpha^\phi_t(1) = \rho^\phi_t(1)$;

\textup{(2)} $\phi(z)\alpha^\phi_t(A) = \alpha^\phi_{st}(A)\phi(z)$;

\textup{(3)} $\rho^\phi_s(S)\phi(z) = \phi(Sz)$;

\textup{(4)} $\rho^\phi_{st}(S\otimes 1) = \rho^\phi_{st}(1)\rho^\phi_s(S)
= \rho^\phi_s(S)\rho^\phi_{st}(1)$; and

\textup{(5)} $\rho^\phi_{ts}(1\otimes S) = \alpha^\phi_t(\rho^\phi_s(S))$.
\end{lemma}

\begin{remark} We have written $S\otimes 1$ for the operator satisfying
$S\otimes 1(xy) = (Sx)y$ for $x\in E_s$ and $y\in E_t$.  If more precision
is needed we write $S\otimes 1^t$.
Similarly, $1\otimes S(yx) = y(Sx)$.
\end{remark}

\begin{proof}[Proof of Lemma~\ref{lemma:alpharho}]
(1) If $\Bb$ is an orthonormal basis for $E_t$, then both
$\alpha^\phi_t(1)$ and $\rho^\phi_t(1)$ equal $\sum_{u\in\Bb} \phi(u)\phi(u)^*$.

(2) See \cite[Lemma~3.6]{fowrae}.

(3) Take a bounded net $\{K_\lambda\}$ of compact operators which converges
strongly to $S$.  Then $K_\lambda \to S$ $\sigma$-weakly as well,
so (3) follows from \eqref{eq:rho}.

(4) By linearity and $\sigma$-weak continuity we may assume that $S=\rankone xy$.
The series $\sum_{u\in\Bb} \rankone{xu}{yu}$
converges $\sigma$-weakly to $S\otimes 1$, so
\begin{align*}
\rho^\phi_{st}(S\otimes 1)
& = \sum_{u\in\Bb} \phi(xu)\phi(yu)^*
  = \phi(x)\Bigl( \sum_{u\in\Bb} \phi(u)\phi(u)^* \Bigr) \phi(y)^* \\
& = \phi(x)\alpha^\phi_t(1)\phi(y)^*
  = \alpha^\phi_{st}(1)\phi(x)\phi(y)^*
  = \rho^\phi_{st}(1)\rho^\phi(S).
\end{align*}
The proof that $\rho^\phi_{st}(S\otimes 1) = \rho^\phi_s(S)\rho^\phi_{st}(1)$ 
is similar.

(5) Again we may assume that $S = \rankone xy$, and
\[
\alpha^\phi_t(\rho^\phi_s(S)) = \sum_{u\in\Bb} \phi(u)\phi(x)\phi(y)^*\phi(u)^*
= \rho^\phi_{ts} \Bigl( \sum_{u\in\Bb} \rankone{ux}{uy} \Bigr) = \rho^\phi_{ts}(1\otimes S).
\]
\end{proof}

We are primarily interested in the case where $P$ is embedded in a group $G$,
and $(G,P)$ is {\em quasi-lattice ordered\/} in the sense of Nica \cite{nica};
see also \cite[\S1]{lacarae} and \cite[\S3]{fowrae}.
Briefly, this means that $P\cap P^{-1} = \{e\}$,
so that $s\le t$ iff $s^{-1}t\in P$ defines a (left-invariant) partial order on $G$,
and that in this partial order every finite subset
of $G$ which has an upper bound in $P$ has a least upper bound in $P$.
When $s,t\in P$ have a common upper bound we write $s\vee t$ for the least upper bound;
when $s$ and $t$ have no common upper bound we write $s\vee t = \infty$.
Our main examples will be direct sums and free products of totally-ordered amenable groups. 

If $(G,P)$ is a quasi-lattice ordered group and
$E$ is a product system over $P$, we say that a
representation $\phi$ of $E$ on $\Hh$ is {\em covariant\/} if
\begin{equation}\label{eq:defn covariance}
\rho^\phi_s(1)\rho^\phi_t(1)
= \begin{cases} \rho^\phi_{s\vee t}(1)
& \text{if $s\vee t < \infty$} \\
0 & \text{otherwise.} \end{cases}
\end{equation}
(By Lemma~\ref{lemma:alpharho}(1), this definition agrees with \cite[Definition~3.2]{fowrae}.)
By \cite[Theorem~4.3]{fowrae}, there is a pair $(\cov(P,E), i_E)$ consisting
of a unital $C^*$-algebra $\cov(P,E)$ and a representation $i_E:E\to\cov(P,E)$
which is universal for covariant representations of $E$, in the sense that

(a) there is a faithful unital representation $\sigma$ of $\cov(P,E)$ on
Hilbert space such that $\sigma\circ i_E$ is a covariant
representation of $E$;

(b) for every covariant representation $\phi$ of $E$ there is a
unital representation $\phi_*$ of $\cov(P,E)$,
called the {\em integrated form of $\phi$\/},
such that $\phi=\phi_*\circ i_E$; and

(c) $\cov(P,E)$ is generated as a $C^*$-algebra by $i_E(E)$.

\noindent Condition (a) falls short of assuring that every representation of $\cov(P,E)$
is the integrated form of a covariant representation.  The following example
shows that there are product systems for which this deficiency is unavoidable.

\begin{example}\label{example}
Let $\Ee$ be a Hilbert space with basis $\{e_k: k\in\NN\}$,
and let $E$ be the product system over $\NN\oplus\NN$ whose fiber over $(m,n)$ is
$\{(m,n)\} \times \Ee^{\otimes(m+n)} \cong \Ee^{\otimes(m+n)}$,
and with multiplication given by
\[
((k,l),\xi)((m,n),\eta) := ((k+m,l+n),\xi\otimes\eta).
\]
If $\{S_k: k\in\NN\}\subset\Bb(\Hh)$ is a collection of isometries with mutually orthogonal ranges,
there is a representation $\phi^S$ of $E$ which satisfies
\[
\phi^S((m,n),e_{k_1}\otimes e_{k_2} \otimes\dotsm\otimes e_{k_{m+n}})
:= S_{k_1}S_{k_2}\dotsm S_{k_{m+n}}.
\]
The projection $L_{m,n}$ onto the closed linear span of $\phi^S(E_{m,n})\Hh$ depends only on $m+n$.
Moreover, the function $m+n\mapsto L_{m,n}$ is strictly decreasing
if $\sum S_kS_k^* < 1$, and is constant if $\sum S_kS_k^* = 1$.
Hence $\phi^S$ is covariant if and only if $\sum S_kS_k^* = 1$.

Suppose then that $\sum S_kS_k^* = 1$.  Since $\phi^S$ is covariant, there is a representation
$\phi^S_*:\cov(P,E)\to\Bb(\Hh)$ such that $\phi^S_*\circ i_E = \phi^S$.
Let $\{T_k: k\in\NN\}$ be another Toeplitz-Cuntz family with $\sum T_kT_k^* < 1$.
By Cuntz's uniqueness theorem for $\Oo_\infty$,
there is a representation $\sigma$ of $\phi^S_*(\cov(P,E)) = C^*(S_k)$
such that $\sigma(S_k) = T_k$.  By the previous paragraph,
$(\sigma\circ\phi^S_*)\circ i_E = \sigma\circ\phi^S = \phi^T$ is not covariant,
and hence $\sigma\circ\phi^S_*$ is not the integrated form of a covariant representation.
\end{example}

The following Proposition hints at how to avoid this kind of pathology.

\begin{prop}\label{prop:covariance}
Suppose $(G,P)$ is a quasi-lattice ordered group and
$E$ is a product system over $P$.
A representation $\phi:E\to\Bb(\Hh)$ is covariant if and only if,
whenever $s,t\in P$, $S\in\Kk(E_s)$, and $T\in\Kk(E_t)$, we have
\begin{equation}\label{eq:covariance}
\rho^\phi_s(S)\rho^\phi_t(T)
= \begin{cases}
\rho^\phi_{s \vee t}((S \otimes 1)(T \otimes 1))
& \text{if $s \vee t < \infty$} \\
0 & \text{if $s \vee t = \infty$.} \end{cases}
\end{equation}
\end{prop}

\begin{proof}
If $\phi$ is covariant, then
\begin{align*}
\rho^\phi_s(S)\rho^\phi_t(T)
& = \rho^\phi_s(S)\rho^\phi_s(1)\rho^\phi_t(1)\rho^\phi_t(T) \\
& = \rho^\phi_s(S)\rho^\phi_{s\vee t}(1)\rho^\phi_t(T)
  = \rho^\phi_{s\vee t}((S\otimes 1)(T\otimes 1)),
\end{align*}
where the last equality uses Lemma~\ref{lemma:alpharho}(4).
The converse follows by $\sigma$-weak continuity.
\end{proof}

If the product $(S\otimes 1)(T\otimes 1)$ in \eqref{eq:covariance}
were known to be compact, this characterization of covariance would make
sense when $\phi$ maps into an abstract $C^*$-algebra.
This motivates the following Definition.

\begin{definition}\label{defn:CA}
Suppose $(G,P)$ is a quasi-lattice ordered group and
$E$ is a product system over $P$.  We say that $E$ is
{\em compactly aligned} if whenever $s,t\in P$ have a common upper
bound and $S$ and $T$ are compact operators on $E_s$ and $E_t$, respectively,
the operator
$(S\otimes 1)(T\otimes 1)$
is a compact operator on $E_{s \vee t}$.
If $E$ is compactly aligned and $\phi$ is a representation of $E$ in
a $C^*$-algebra $B$, we say that $\phi$ is {\em covariant\/} if
\eqref{eq:covariance} holds
whenever $s,t\in P$, $S\in\Kk(E_s)$ and $T\in\Kk(E_t)$.
\end{definition}

\begin{remark}
If $E$ has finite-dimensional fibers,
or if $(G,P)$ is a total order,
then $E$ is obviously compactly aligned.
\end{remark}

The following Proposition shows that the pathology observed in Example~\ref{example}
cannot occur when $E$ is compactly aligned; 
every representation of $\cov(P,E)$ is the integrated form of a covariant representation.

\begin{prop}\label{prop:i_E covariant}
Let $E$ be a compactly-aligned product system over $P$.

\textup{(1)} The universal map $i_E:E\to\cov(P,E)$ is covariant.

\textup{(2)} If $\sigma$ is a representation of $\cov(P,E)$,
then $\sigma\circ i_E$ is covariant.
\end{prop}

\begin{proof} (1) Let $\sigma$ be a faithful representation of $\cov(P,E)$
on a Hilbert space such that $\phi := \sigma\circ i_E$ is covariant.
Since $E$ is compactly aligned we can apply \eqref{eq:compose}
to both sides of \eqref{eq:covariance} to obtain
\begin{equation}\label{eq:iE}
\sigma\bigl(\rho^{i_E}_s(S)\rho^{i_E}_t(T)\bigr)
= \begin{cases}
\sigma\bigl(\rho^{i_E}_{s \vee t}((S \otimes 1)(T \otimes 1))\bigr)
& \text{if $s \vee t < \infty$} \\
0 & \text{if $s \vee t = \infty$.} \end{cases}
\end{equation}
The result follows by applying $\sigma^{-1}$ to this equation.

(2) Taking $\phi = i_E$ in \eqref{eq:covariance} and applying $\sigma$ to both sides
gives \eqref{eq:iE}; covariance of $\sigma\circ i_E$ then follows from \eqref{eq:compose}.
\end{proof}

For compactly-aligned product systems, we can improve on \cite[Proposition~3.7]{fowrae};
the series given there converges in {\em norm\/}:

\begin{prop}\label{prop:v*w} Let $E$ be a compactly-aligned product system over $P$,
and suppose $v,w\in E$ satisfy $p(v)\vee p(w) < \infty$.
Let $\Bb$ and $\Cc$ be orthonormal bases for $E_{p(v)^{-1}(p(v)\vee p(w))}$
and $E_{p(w)^{-1}(p(v)\vee p(w))}$, respectively.

\textup{(1)} The series
\[
\sum_{f\in\Bb, g\in\Cc} \langle wg,vf \rangle \rankone{vf}{wg}
\]
converges in norm to
$((\rankone vv) \otimes 1)((\rankone ww) \otimes 1)
\in\Kk(E_{p(v)\vee p(w)})$.

\textup{(2)} If $\phi:E\to B$ is covariant,
then the series
\[
\sum_{f\in\Bb, g\in\Cc} \langle wg,vf \rangle \phi(f)\phi(g)^*
\]
converges in norm to $\phi(v)^*\phi(w)$.
\end{prop}

\begin{proof} Since $E$ is compactly aligned,
$K := ((\rankone vv)\otimes 1)((\rankone ww)\otimes 1)$
is compact.
The finite sums $\sum 1\otimes(\rankone ff)$ are finite-rank projections which
increase to $1$, so the series
$\sum_{f\in\Bb} (1\otimes(\rankone ff))K$ converges in norm to $K$.
Similarly, the series $\sum_{g\in\Cc} K(1\otimes(\rankone gg))$ converges
in norm to $K$.  Hence
\[
\sum_{f\in\Bb, g\in\Cc} (1\otimes(\rankone ff))K(1\otimes(\rankone gg))
\]
converges in norm to $K$.  Since
\begin{align*}
(1\otimes(\rankone ff))K(1\otimes(\rankone gg))
& = ((\rankone vv)\otimes(\rankone ff))((\rankone ww)\otimes(\rankone gg)) \\
& = (\rankone{vf}{vf})(\rankone{wg}{wg}) \\
& = \langle wg, vf \rangle \rankone{vf}{wg},
\end{align*}
this gives (1).

(2) The series
\[
\sum_{f\in\Bb, g\in\Cc} \langle wg,vf \rangle \phi(vf)\phi(wg)^*
=
\sum_{f\in\Bb, g\in\Cc} \langle wg,vf \rangle \rho^\phi_{p(v)\vee p(w)}(\rankone{vf}{wg})
\]
converges in norm to
$\rho^\phi_{p(v)\vee p(w)}(K)$, which by the covariance of $\phi$ is equal to
\[
\rho^\phi_{p(v)}(\rankone vv)\rho^\phi_{p(w)}(\rankone ww)
= \phi(v)\phi(v)^*\phi(w)\phi(w)^*.
\]
Multiplying on the left by $\phi(v)^*$ and on the right by $\phi(w)$ gives (2).
\end{proof}

The analysis of many $C^*$-algebras is facilitated by the existence of a
set of Wick-ordered monomials which have dense linear span.
The following Corollary shows that when $E$ is compactly aligned,
$\cov(P,E)$ is such an algebra.

\begin{cor}\label{cor:wick}
If $E$ is a compactly-aligned product system over $P$, then
\[
\cov(P,E) = \clsp\{i_E(x)i_E(y)^*: x,y\in E\}.
\]
\end{cor}

\begin{proof} Taking $\phi = i_E$ in Proposition~\ref{prop:v*w}(2) shows that
the right hand side is closed under multiplication.
The result follows since each generator $i_E(x)$ of $\cov(P,E)$ can be written as
$i_E(x) = i_E(x)i_E(\Omega)^*$, where $\Omega$ is the identity of $E$.
\end{proof}

Since our results in the next section involve free products,
we need the following Proposition.
(See \cite[Examples~1.4(E4)]{fowrae} for the definition of the free product.)

\begin{prop}\label{prop:free products}
For each $\lambda$ belonging to some index set $\Lambda$,
let $(G^\lambda, P^\lambda)$ be a quasi-lattice ordered group,
and let $E^\lambda$ be a  product system over $P^\lambda$.
If each $E^\lambda$ is compactly aligned, then so is their free product $*E^\lambda$.
\end{prop}

\begin{proof} Suppose $s,t\in*P^\lambda$ have a common upper bound, $S\in\Kk(E_s)$,
and $T\in\Kk(E_t)$.  If $s \le t$, then the product $(S\otimes 1)(T\otimes 1)$
reduces to $(S\otimes 1)T$, and is hence compact; the case $t\le s$ is similar.
If $s,t < s\vee t$, then there exists $\mu\in\Lambda$ such that $s^{-1}t\in G^\mu$;
that is, we can write $s = rs'$, $t = rt'$ with $r\in*P^\lambda$,
$s',t'\in P^\mu$ and $s'\vee t' < \infty$.
By linearity and continuity we may assume that $S = R_1\otimes S'$,
where $R_1\in\Kk(E_r)$ and $S'\in\Kk(E_{s'})$; similarly, we assume that $T = R_2\otimes T'$.
Then
\begin{align*}
(S\otimes 1)(T\otimes 1)
& = ((R_1\otimes S')\otimes 1)((R_2\otimes T')\otimes 1) \\
& = R_1R_2 \otimes((S'\otimes 1)(T' \otimes 1))
\end{align*}
is compact since $E^\mu$ is compactly aligned.
\end{proof}

\begin{remark}
In \cite[Examples~1.4]{fowrae}, it was also shown how one can twist a product system $E$
by a multiplier $\omega:P\times P\to\TT$;
it is easy to see that $E^\omega$ is compactly aligned if $E$ is.
One can also take tensor products of product systems,
and this also preserves the property of being compactly aligned.
\end{remark}

We close this section by showing how to adapt Dinh's construction of von Neumann
discrete product systems \cite[Example~6.1]{dinhjfa}
to an arbitrary left-cancellative semigroup $P$.
Fix a separable Hilbert space $\Hh$ and a unit vector $\xi\in\Hh$.
For $s\in P$, let $\Hh_s$ be the von Neumann tensor product
$\otimes_{P\setminus sP} \Hh$ with canonical unit vector $\xi$;
that is, $\Hh_s$ is the Hilbert space inductive limit $\varinjlim \otimes_F \Hh$,
taken over the finite subsets $F\subseteq P\setminus sP$, under the isometric
embeddings obtained by tensoring with $\xi$ \cite{vonNeumann}.
Vectors of the form $\otimes_{r\in P\setminus sP} \eta_r$, where $\eta_r = \xi$
for all but finitely many $r$, have dense linear span in $\Hh_s$.

Let $E := \bigsqcup_{s\in P} \{s\}\times\Hh_s$, $p(s,\eta) := s$,
and define multiplication in $E$ by
\[
\Bigl( s, \bigotimes_{r\in P\setminus sP} \eta_r \Bigr)
\Bigl( t, \bigotimes_{r\in P\setminus tP} \zeta_r \Bigr)
:=
\Bigl( st, \bigotimes_{r\in P\setminus stP} \gamma_r \Bigr),
\]
where
\[
\gamma_r := \begin{cases} \eta_r & \text{if $r\in P \setminus sP$,} \\
\zeta_a & \text{if $r = sa$ for some $a\in P\setminus tP$.} \end{cases}
\]
It is routine to check that $p:E\to P$ is a product system.
Let $d$ be the dimension of $\Hh$. 
Since any other pair $(\Hh',\xi')$ with $\dim\Hh' = d$ determines
a product system which is isomorphic to the one above, we denote
this product system $E^d$, and call it a {\em von Neumann discrete product system.}

\begin{prop}\label{prop:von Neumann} If $(G,P)$ is a quasi-lattice ordered group,
then any von Neumann discrete product system over $P$ is compactly aligned.
\end{prop}

\begin{proof} Suppose $S\in\Kk(E_s)$, $T\in\Kk(E_t)$, and $s\vee t < \infty$;
we will show that $(S\otimes 1)(T\otimes 1)$ is compact.
Since $P \setminus sP = (P\setminus(sP\cup tP)) \sqcup (tP \setminus sP)$,
there is a canonical isomorphism
\[
E_s \cong \Bigl( \bigotimes_{P \setminus (sP\cup tP)} \Hh \Bigr)
          \otimes
          \Bigl( \bigotimes_{tP \setminus sP}\Hh \Bigr).
\]
Similarly, there is a canonical isomorphism
\[
E_t \cong \Bigl( \bigotimes_{P \setminus (sP\cup tP)} \Hh \Bigr)
          \otimes
          \Bigl( \bigotimes_{sP \setminus tP}\Hh \Bigr).
\]
By linearity and continuity, we may assume that under these identifications
$S = S_1 \otimes S_2$ and $T = T_1 \otimes T_2$, where $S_1$, $S_2$, $T_1$, and $T_2$
are compact.  Using the decomposition
$P \setminus (s\vee t)P
   = (P\setminus(sP\cup tP)) \sqcup (tP \setminus sP) \sqcup (sP \setminus tP)$,
we find that $S\otimes 1 = S_1\otimes S_2 \otimes 1$
and $T\otimes 1 = T_1 \otimes 1 \otimes T_2$,
so
\[
(S\otimes 1)(T\otimes 1) = S_1T_1 \otimes S_2 \otimes T_2
\]
is compact.
\end{proof}

\section{Faithful representations}\label{section:faithful}

Proposition~\ref{prop:v*w} allows us to apply the main results of \cite{fowrae}
to compactly-aligned product systems.  For this we need a technical amenability hypothesis
on $(G,P)$; namely, we require the existence of a homomorphism $\theta:(G,P)\to(\Gg,\Pp)$ of
quasi-lattice ordered groups such that $\Gg$ is amenable and,
whenever $s\vee t < \infty$, we have
\begin{equation}\label{eq:amenable}
\theta(s\vee t) = \theta(s) \vee \theta(t)
\quad\text{and}\quad
\theta(s) = \theta(t) \Rightarrow s = t.
\end{equation}
This hypothesis is satisfied trivially if $G$ is an amenable group,
and if $(G,P)$ is a free product of quasi-lattice ordered groups $(G^\lambda, P^\lambda)$
with each $G^\lambda$ amenable, we can take $\theta$ to be the canonical map
$*(G^\lambda,P^\lambda)\to\oplus(G^\lambda,P^\lambda)$
\cite[Proposition~4.3]{lacarae}.

\begin{theorem}\label{theorem:sufficient}
Suppose $(G,P)$ is a quasi-lattice ordered group
which admits a homomorphism $\theta:(G,P)\to(\Gg,\Pp)$ as above.
If $E$ is a compactly-aligned product system over $P$
and $\phi$ is a covariant representation of $E$ on Hilbert space such that
\begin{equation}\label{eq:sufficient}
\prod_{k=1}^n (1 - \rho^\phi_{s_k}(1)) \ne 0
\quad\text{whenever $s_1, \dots, s_n\in P\setminus\{e\}$,}
\end{equation}
then the integrated form $\phi_*$ is faithful on $\cov(P,E)$.
\end{theorem}

\begin{proof} Since $E$ is compactly aligned, Proposition~\ref{prop:v*w}
implies that the hypotheses of \cite[Theorem~6.1]{fowrae} are satisfied.
This in turn implies that the hypotheses of \cite[Theorem~5.1]{fowrae} are satisfied,
so \eqref{eq:sufficient} implies that $\pi_\phi\times\phi$ is a faithful representation of $\bpp$.
Since $\phi_*$ is the restriction of $\pi_\phi\times\phi$
to the subalgebra $\cov(P,E)$ (see \cite[Theorem~4.3]{fowrae}),
it is faithful when \eqref{eq:sufficient} holds.
\end{proof}

The following Corollary generalizes \cite[Corollaries~5.2 and 5.3]{lacarae}.
The second part involves quasi-lattice orders $(G,P)$ which satisfy the following
condition:
\begin{multline}\label{eq:div}
\text{for each finite subset $F\subseteq P\setminus\{e\}$,} \\
\text{there exists $s\in P\setminus\{e\}$ such that $s < s\vee t$ for every $t\in F$.}
\end{multline}
For example, one could take $G$ to be a countable dense subgroup of $\RR$ with positive cone
$P = G\cap [0,\infty)$; one could also take a direct sum of such $(G,P)$.

\begin{cor}\label{cor:simplicity}
For each $\lambda$ belonging to some index set $\Lambda$,
let $(G^\lambda,P^\lambda)$ be a quasi-lattice ordered group with $G^\lambda$ amenable,
and let $E^\lambda$ be a compactly-aligned product system over $P^\lambda$.
Then $\cov(*P^\lambda, *E^\lambda)$ is simple if

\textup{(1)} $\Lambda$ is infinite, or if

\textup{(2)} $\abs\Lambda \ge 2$,
and there exists $\mu\in\Lambda$ such that $(G^\mu,P^\mu)$
satisfies \eqref{eq:div}.
\end{cor}

\begin{proof} By Proposition~\ref{prop:free products}, $*E^\lambda$ is compactly aligned,
so by Proposition~\ref{prop:i_E covariant}, every representation of $\cov(*P^\lambda,*E^\lambda)$
is the integrated form of a covariant representation.  Since the canonical homomorphism
$\theta:*(G^\lambda,P^\lambda) \to \oplus (G^\lambda,P^\lambda)$ satisfies \eqref{eq:amenable},
Theorem~\ref{theorem:sufficient} applies, and it suffices to verify that \eqref{eq:sufficient}
holds for every covariant representation $\phi$.

Suppose $s_1, \dots, s_n\in *P^\lambda\setminus\{e\}$.
Express $s_k = r_ks_k'$ with $s_k'\in *P^\lambda$
and $r_k \in P^{\lambda_k}\setminus\{e\}$ for some $\lambda_k\in\Lambda$.
Since $\rho^\phi_{s_k}(1) \le \rho^\phi_{r_k}(1)$, it suffices to show that
\begin{equation}\label{eq:sufficient2}
\prod_{k=1}^n (1 - \rho^\phi_{r_k}(1)) \ne 0.
\end{equation}

Suppose (1) holds.  Then there exists $\lambda\in\Lambda\setminus\{\lambda_k\}$,
and for any $t\in P^\lambda\setminus\{e\}$ we have $t\vee r_k = \infty$ for each $k$.
Since $\phi$ is covariant we deduce that
$\rho^\phi_t(1) \prod_{k=1}^n (1 - \rho^\phi_{r_k}(1))
= \rho^\phi_t(1) \ne 0$,
and this implies \eqref{eq:sufficient2}.

Next suppose (2) holds.
By reordering if necessary, assume that there exists $m \le n$ such that
$r_k \in P^\mu$ if and only if $1 \le k \le m$.
By \eqref{eq:div}, there exists $r\in P^\mu\setminus\{e\}$ such that
$r < r\vee r_k$ for $1 \le k \le m$; observe that $r\vee r_k = \infty$
if $k > m$.  By the covariance of $\phi$ and Lemma~\ref{lemma:alpharho},
\[
\rho^\phi_r(1) \prod_{k=1}^n (1 - \rho^\phi_{r_k}(1))
 = \alpha^\phi_r\Bigl( \prod_{k=1}^m (1 - \rho^\phi_{r^{-1}(r\vee r_k)}(1)) \Bigr).
\]
Now $\abs\Lambda \ge 2$, so there exists $\lambda\in\Lambda\setminus\{\mu\}$.
Let $t\in P^\lambda\setminus\{e\}$.
Then $t\vee s = \infty$ for every $s\in P^\mu\setminus\{e\}$,
so multiplying both sides of the previous equation
by $\alpha^\phi_r(\rho^\phi_t(1))$ gives
\begin{align*}
\alpha^\phi_r(\rho^\phi_t(1)) \prod_{k=1}^n (1 - \rho^\phi_{r_k}(1))
& = \alpha^\phi_r\Bigl( \rho^\phi_t(1) \prod_{k=1}^m (1 - \rho^\phi_{r^{-1}(r\vee r_k)}(1)) \Bigr) \\
& = \alpha^\phi_r(\rho^\phi_t(1)) \ne 0.
\end{align*}
This implies \eqref{eq:sufficient2}, completing the proof.
\end{proof}

When $E$ has infinite-dimensional fibers, \eqref{eq:sufficient} may not be necessary
for faithfulness.  For example, if $E$ is the product system over $\NN$ such that
$\dim E_1 = \infty$, then $\cov(\NN,E)$ is the simple $C^*$-algebra $\Oo_\infty$,
whereas the previous Theorem gives only that a Toeplitz-Cuntz family $\{S_1, S_2, \dots\}$
generates a faithful copy of $\cov(\NN,E)$ if $\sum S_kS_k^* < 1$.
(See \cite[Corollary~1.6 and Examples~5.6(2)]{fowrae} for details.)

To obtain a theorem which completely characterizes the faithful representations
of $\cov(P,E)$ when $E$ has infinite-dimensional fibers, we assume that
$(G,P)$ is a free product of totally-ordered amenable groups.
Then
\begin{equation}\label{eq:qto}
s\vee t\in\{s,t,\infty\}\qquad\text{for every $s,t\in P$,}
\end{equation}
a property we shall use extensively.
Notice that \eqref{eq:qto} implies that any product system over $P$ is compactly aligned:
if $S\in\Kk(E_s)$ and $T\in\Kk(E_t)$ with $s\vee t < \infty$, then the product
$(S\otimes 1)(T\otimes 1)$ reduces to either $S(T\otimes 1)$ or $(S\otimes 1)T$,
and is hence compact.

\begin{theorem}\label{theorem:faithful}
For each $\lambda$ belonging to some index set $\Lambda$,
let $(G^\lambda,P^\lambda)$ be a totally-ordered amenable group.
Let $E$ be a product system over $P:= *P^\lambda$,
and let $\phi$ be a covariant representation of $E$ in a unital $C^*$-algebra $B$.
Then $\phi_*:\cov(P,E)\to B$ is injective if and only if
\begin{equation}\label{eq:faithful}
\begin{split}
\prod_{k=1}^n (1 - \rho^\phi_{s_k}(1)) \ne 0
\quad\text{whenever $\{s_1,\dots,s_n\}$ is a finite subset} \\
\text{of $P\setminus\{e\}$ and $\dim E_{s_k} < \infty$ for each $k$.}
\end{split}
\end{equation}
\end{theorem}

\begin{remark}
Our proof does not use the full strength of these hypotheses,
and so may apply in slightly more generality:
we require only that $(G,P)$ satisfies \eqref{eq:qto},
and that it admits a homomorphism $\theta:(G,P)\to(\Gg,\Pp)$
of the type discussed prior to Theorem~\ref{theorem:sufficient}.
\end{remark}

One direction of this theorem is easy: \eqref{eq:faithful} is satisfied
by the left regular representation $l:E\to\Bb(\bigoplus_{s\in P} E_s)$
\cite[Lemmas~1.10 and~3.5]{fowrae},
and since the identity operator on each $E_{s_k}$ is compact,
we may apply \eqref{eq:compose} to obtain
$1 - \rho^l_{s_k}(1) = l_*(1 - \rho^{i_E}_{s_k}(1))$;
thus \eqref{eq:faithful} is satisfied by the universal map $i_E$.
Composing with an injective $\phi_*$ and again using \eqref{eq:compose}
gives \eqref{eq:faithful}.

For the converse, we model our proof after \cite[Theorem~5.1]{fowrae},
making use of the identification of $\cov(P,E)$ as a subalgebra of $\bpp$
\cite[Theorem~4.3]{fowrae}.
Proposition~\ref{prop:v*w} and the canonical map
$\theta:*(G^\lambda,P^\lambda)\to\oplus(G^\lambda,P^\lambda)$
allow us to conclude from \cite[Theorem~6.1]{fowrae} that $E$
is amenable; this means that the canonical conditional expectation
$\Phi_\delta$ on $\bpp$ is faithful on positive elements \cite[p.~189]{fowrae}.
Since this expectation satisfies
\begin{equation}\label{eq:ce}
\Phi_\delta(i_E(x)i_E(y)^*)
= \begin{cases} i_E(x)i_E(y)^* & \text{if $p(x) = p(y)$} \\
0 & \text{otherwise,}
\end{cases}
\end{equation}
Corollary~\ref{cor:wick} shows that
$\Phi_\delta$ restricts to a map on $\cov(P,E)$;
indeed, \eqref{eq:ce} determines $\Phi_\delta$ on $\cov(P,E)$.

In Proposition~\ref{prop:fpa} we show that $\phi_*$ is isometric on
$\Phi_\delta(\cov(P,E))$,
and we construct a spatial version $\Phi_\phi$
of $\Phi_\delta$ such that $\Phi_\phi\circ\phi_* = \phi_*\circ\Phi_\delta$.
A standard argument then shows that $\phi_*$ is injective: $\phi_*(b) = 0
\Rightarrow \Phi_\phi\circ\phi_*(b^*b) = 0
\Rightarrow \phi_*\circ\Phi_\delta(b^*b) = 0
\Rightarrow \Phi_\delta(b^*b) = 0
\Rightarrow b = 0$.

For the remainder of this paper we will assume that the hypotheses
of Theorem~\ref{theorem:faithful} are satisfied, and that \eqref{eq:faithful} holds.
We also assume that $B$ is represented
faithfully and nondegenerately on a Hilbert space $\Hh$,
and thus regard $\phi$ as a covariant representation of $E$ on $\Hh$.
This allows us to make use of the endomorphisms $\alpha^\phi_s$
defined in \eqref{eq:alpha},
and to extend the domain of $\rho^\phi_s$ to $\Bb(E_s)$.

The following two technical lemmas will be needed to prove Proposition~\ref{prop:fpa}.
The first provides a projection which can be used to kill those terms
of a finite sum $\sum \phi(x_j)\phi(y_j)^*$ which are ``too long'',
and the second provides a vector which plays the role of the aperiodic sequence
used by Cuntz to kill off-diagonal terms \cite{cun}.

\begin{lemma}\label{lemma:technical}
Let $a\in P$,
and suppose $F$ is a finite subset of $E$ such that $p(x)\not\le a$ for every $x\in F$.

\textup{(1)} For every $\epsilon > 0$, there is a nonzero projection $Q_\epsilon\in\Bb(\Hh)$
such that
\[
\norm{\alpha^\phi_a(Q_\epsilon)\phi(x)} < \epsilon
\qquad\text{for every $x\in F$.}
\]

\textup{(2)} More specifically, suppose $C$ is a finite subset of $P\setminus\{e\}$
such that $a^{-1}p(x)\in C$ for every $x\in F$ such that $a < p(x)$.
Let $R_C$ be the projection $\prod_{r\in C} (1 - \rho^\phi_r(1))$.
Then $\alpha^\rho_a(R_C)\phi(x) = 0$ for every $x\in F$, and if $R_C = 0$,
then there is a unit vector $y_\epsilon\in E$ such that
\[
\norm{\alpha^\phi_a(\phi(y_\epsilon)\phi(y_\epsilon)^*)\phi(x)} < \epsilon
\qquad\text{for every $x\in F$.}
\]
\end{lemma}

\begin{proof} We need only prove (2).
If $x\in F$ satisfies $a \vee p(x) = \infty$,
then for any $Q\in\Bb(\Hh)$,
Lemma~\ref{lemma:alpharho} and the covariance of $\phi$ give
\[
\alpha^\phi_a(Q)\phi(x) = \alpha^\phi_a(Q)\rho^\phi_a(1)\rho^\phi_{p(x)}(1)\phi(x) = 0.
\]
Hence we can assume that $a \vee p(x) < \infty$, and thus $a < p(x)$, for every $x\in F$.

If $x\in F$, then $r := a^{-1}p(x)\in C$,
and for any $R\in\Bb(E_r)$ we can use Lemma~\ref{lemma:alpharho} to calculate
\begin{equation}\label{eq:norm}
\begin{split}
\norm{\alpha^\phi_a(1 - \rho^\phi_r(R))\phi(x)}
& = \norm{(\rho^\phi_a(1) - \rho^\phi_{ar}(1^a\otimes R))\rho^\phi_{ar}(1)\phi(x)} \\
& = \norm{\rho^\phi_{ar}(1 - (1^a\otimes R))\phi(x)} \\
& = \norm{\phi((1 - (1^a\otimes R))x)} \\
& = \norm{(1^a \otimes (1 - R))x};
\end{split}
\end{equation}
in particular $\alpha^\phi_a(1 - \rho^\phi_r(1))\phi(x) = 0$.
Thus $\alpha^\phi_a(R_C)\phi(x) = 0$ for every $x\in F$.

Suppose $R_C = 0$.
Let $C_{\min}$ be the set of minimal elements of $C$.
Since the projections $1 - \rho^\phi_r(1)$ increase in $r$, we have
\[
0 = \prod_{r\in C} (1 - \rho^\phi_r(1)) = \prod_{r\in C_{\min}} (1 - \rho^\phi_r(1)).
\]
By \eqref{eq:faithful}, there exists $r_{\min}\in C_{\min}$ such that
$\dim E_{r_{\min}} = \infty$.
Let $r_{\max}$ be a maximal element of $\{r\in C: r_{\min} \le r\}$.
If $r\in C$ and $r\not\le r_{\max}$, then by \eqref{eq:qto} and the maximality of
$r_{\max}$ we have $r\vee r_{\max} = \infty$;
for any $y\in E_{r_{\max}}$, the covariance of $\phi$ thus gives
\[
\phi(y)^*\rho^\phi_r(1)
= \phi(y)^*\rho^\phi_{r_{\max}}(1)\rho^\phi_r(1)
= 0.
\]
Hence $\phi(y)\phi(y)^* \le 1 - \rho^\phi_r(1)$,
and from \eqref{eq:norm} we deduce that
\begin{equation}\label{eq:rmax}
\alpha^\phi_a(\phi(y)\phi(y)^*)\phi(x) = 0
\qquad\text{if $y\in E_{r_{\max}}$, $x\in F$, and $a^{-1}p(x) \not\le r_{\max}$.}
\end{equation}

By \eqref{eq:qto}, the set $\{r\in C: r \le r_{\max}\}$
is totally ordered, say $r_{\min} = r_1 < r_2 < \dots < r_{n-1} < r_n = r_{\max}$.
For each $k\in\{1,\dots, n\}$, let $P_k$ be a finite-rank projection
on $E_{r_k}$ which satisfies
\[
\norm{(1^a\otimes(1 - P_k))x} < \epsilon
\qquad\text{for every $x\in F$ such that $p(x) = ar_k$.}
\]
By \eqref{eq:norm} and \eqref{eq:rmax}, it suffices to find
a unit vector $y_\epsilon\in E_{r_{\max}}$ such that
$\phi(y_\epsilon)\phi(y_\epsilon)^* \le 1 - \rho^\phi_{r_k}(P_k)$
for every $k\in\{1,\dots, n\}$; that is, such that
\begin{equation}\label{eq:yepsilon}
\phi(y_\epsilon)^*\rho^\phi_{r_k}(P_k) = 0
\qquad\text{for every $k\in\{1,\dots, n\}$.}
\end{equation}

If some $E_{r_k^{-1}r_{k+1}}$ is finite-dimensional, then we can create a new
collection of finite-rank projections by removing $r_k$ and replacing
$P_{k+1}$ by $P_{k+1} \vee (P_{r_k}\otimes 1^{r_k^{-1}r_{k+1}})$.
For any $y\in E_{r_{\max}}$ we have $\phi(y)^* = \phi(y)^*\rho^\phi_{r_{k+1}}(1)$,
so Lemma~\ref{lemma:alpharho}(4) gives
\[
  \phi(y)^* \rho^\phi_{r_k}(P_k)
= \phi(y)^* \rho^\phi_{r_{k+1}}(1) \rho^\phi_{r_k}(P_k)
= \phi(y)^* \rho^\phi_{r_{k+1}}(P_{r_k}\otimes 1^{r_k^{-1}r_{k+1}});
\]
hence any $y_\epsilon$ which satisfies \eqref{eq:yepsilon}
for this new collection also works for the original.
We therefore assume that each $E_{r_k^{-1}r_{k+1}}$ is infinite-dimensional.

We will define unit vectors $w_1, \dots, w_n\in E$ recursively such that
\begin{equation}\label{eq:recursion}
w_k\in E_{r_k},
\quad\text{ and }\quad
\phi(w_k)^*\rho^\phi_{r_j}(P_{r_j}) = 0
\quad\text{ for $j\in\{1,\dots, k\}$;}
\end{equation}
then $y_\epsilon:=w_n$ satisfies \eqref{eq:yepsilon}, completing the proof.
Since $\dim E_{r_1} = \infty$, we can choose $w_1\in E_{r_1}$
such that $P_{r_1}w_1 = 0$; then
\[
\phi(w_1)^*\rho^\phi_{r_1}(P_{r_1})
 = \phi(w_1)^*\rho^\phi_{r_1}((\rankone{w_1}{w_1})P_{r_1}) = 0.
\]
Suppose $1 \le k \le n-1$ and $w_k$ satisfies \eqref{eq:recursion}.
We claim that there is a unit vector $z\in E_{r_k^{-1}r_{k+1}}$ such that
$P_{k+1}(w_kz) = 0$.  Given this, $w_{k+1} := w_kz$ satisfies \eqref{eq:recursion},
since
\[
\phi(w_kz)^*\rho^\phi_{r_{k+1}}(P_{r_{k+1}})
= \phi(w_kz)^*\rho^\phi_{r_{k+1}}((\rankone{w_kz}{w_kz})P_{r_{k+1}}) = 0,
\]
and, for $1 \le j \le k$,
\[
\phi(w_kz)^*\rho^\phi_{r_j}(P_{r_j}) = \phi(z)^*\phi(w_k)^*\rho^\phi_{r_j}(P_{r_j}) = 0.
\]

For the existence of such a $z$, let $e_1, \dots, e_m$ be an orthonormal basis
for the range of $P_{k+1}$, so that $P_{k+1} = \sum_{l=1}^m \rankone{e_l}{e_l}$.
For each $l\in\{1,\dots, m\}$,
let $z_l$ be the unique vector in $E_{r_k^{-1}r_{k+1}}$ such that
$((\rankone{w_k}{w_k})\otimes 1^{r_k^{-1}r_{k+1}})e_l = w_kz_l$.
Since $E_{r_k^{-1}r_{k+1}}$ is infinite-dimensional, there is a unit vector
$z\in E_{r_k^{-1}r_{k+1}}$ such that $\langle z, z_l \rangle = 0$ for each $l$.
Then
\begin{align*}
P_{k+1}(w_kz)
& = \sum_{l=1}^m (\rankone{e_l}{e_l})((\rankone{w_k}{w_k})\otimes 1^{r_k^{-1}r_{k+1}})(w_kz) \\
& = \sum_{l=1}^m (\rankone{e_l}{w_kz_l})(w_kz) = 0.
\end{align*}
\end{proof}

\begin{lemma}\label{lemma:aperiodic}
Let $C$ and $D$ be finite subsets of $P\setminus\{e\}$.
If the projection $\prod_{r\in C} (1 - \rho^\phi_r(1))$ vanishes,
then for every unit vector $y\in E$ there is a unit vector $z\in E$
such that
\begin{equation}\label{eq:aperiodic}
\phi(yz)\phi(yz)^*\alpha^\phi_d\bigl(\phi(yz)\phi(yz)^*\bigr) = 0
\qquad\text{for every $d\in D$.}
\end{equation}
\end{lemma}

\begin{proof} First suppose $D$ is a singleton $\{d\}$.
Since
\[
\alpha^\phi_d\bigl(\phi(yz)\phi(yz)^*\bigr) = \sum_{u\in\Bb} \phi(uyz)\phi(uyz)^*,
\]
where $\Bb$ is an orthonormal basis for $E_d$, it suffices to find $z$ such that
\begin{equation}\label{eq:shift}
\phi(yz)^*\phi(uyz) = 0\qquad\text{for every $u\in E_d$.}
\end{equation}
If $p(y) \vee dp(y) = \infty$, then
$\phi(y)^*\phi(uy) = 0$, and \eqref{eq:shift} holds for any vector $z$.
Assume then that $p(y) \vee dp(y) < \infty$, and define $s\in P\setminus\{e\}$ by
\[
s := \begin{cases} p(y)^{-1}dp(y) & \text{if $p(y) < dp(y)$,} \\
                   (dp(y))^{-1}p(y) & \text{if $dp(y) < p(y)$.} \end{cases}
\]
Let $z_1$ be a unit vector in $E_s$.  All we need to achieve \eqref{eq:shift}
is unit vector $z_2\in E$ such that $\phi(z_2)^*\phi(z_1) = 0$;
for then we let $z := z_1z_2$, and calculate
\begin{align*}
\phi(yz)^*\phi(uyz)
& = \phi(z_2)^*\phi(z_1)^*\phi(y)^*\phi(uy)\phi(z_1)\phi(z_2) \\
& = \begin{cases} \langle uy,yz_1 \rangle \phi(z_2)^*\phi(z_1)\phi(z_2)
                  & \text{if $p(y) < dp(y)$,} \\
                  \langle uyz_1,y \rangle \phi(z_2)^*\phi(z_1)^*\phi(z_2)
                  & \text{if $dp(y) < p(y)$} \end{cases} \\
& = 0.
\end{align*}
If there exists $t\in P$ such that $s\vee t = \infty$, then any $z_2\in E_t$
will suffice.
On the other hand, if $s\vee t < \infty$ for every $t\in P$, then
for every $r\in C$ we have either $r \le s$ or $s \le r$, and thus
\[
0 = (1 - \rho^\phi_s(1)) \prod_{r\in C} (1 - \rho^\phi_r(1))
 = 1 - \rho^\phi_t(1)
\]
for some $t \le s$.
The hypothesis \eqref{eq:faithful} implies that $E_t$, and hence $E_s$,
is infinite-dimensional, and taking any $z_2\in E_s$ orthogonal to $z_1$ suffices.

Now suppose inductively that the lemma holds for some finite set $D$,
and that $d'\in P\setminus\{e\}$; we will show that the lemma holds for $D\cup\{d'\}$.
Fix a unit vector $y\in E$, and let $z\in E$ be a unit vector such that
\eqref{eq:aperiodic} holds.
By the above, there exists a unit vector $z'\in E$ such that
\[
\phi(yzz')\phi(yzz')^*\alpha^\phi_{d'}\bigl(\phi(yzz')\phi(yzz')^*\bigr) = 0.
\]
Since $\phi(yzz')\phi(yzz')^* \le \phi(yz)\phi(yz)^*$,
\eqref{eq:aperiodic} holds with $zz'$ in place of $z$, and the induction is complete.
\end{proof}

\begin{prop}\label{prop:fpa}
If $\phi$ satisfies \eqref{eq:faithful}, then

\textup{(1)} $\phi_*$ is isometric on $\Phi_\delta(\cov(P,E))$, and

\textup{(2)} there is a contractive linear map $\Phi_\phi$ on $\phi_*(\cov(P,E))$
such that
\[
\Phi_\phi\circ\phi_* = \phi_*\circ\Phi_\delta.
\]
\end{prop}

\begin{proof} (1) By Corollary~\ref{cor:wick} and the continuity of
$\Phi_\delta$, finite sums of the form
\[
X := \sum_j i_E(x_j)i_E(y_j)^*,
\]
where $p(x_j) = p(y_j)$ for each $j$,
are dense in $\Phi_\delta(\cov(P,E))$.
Hence it suffices to fix such an $X$ and show that $\norm{\phi_*(X)} = \norm X$.

In the proof of \cite[Proposition~5.4]{fowrae}
it was shown that there exists $a\in P$ such that $\norm X = \norm{T_a}$, where
$T_a$ is the operator on $E_a$ defined by
\[
T_a := \sum_{\{j: p(x_j) \le a\}} (\rankone{x_j}{y_j})\otimes 1^{p(x_j)^{-1}a}.
\]
Let $\epsilon > 0$.  By Lemma~\ref{lemma:technical}(1),
there is a nonzero projection $Q_\epsilon\in\Bb(\Hh)$ such that
$\norm{\alpha^\phi_a(Q_\epsilon)\phi(x_j)} < \epsilon$
whenever $p(x_j) \not\le a$.
Since $\phi_*(X) = \sum_j \phi(x_j)\phi(y_j)^*$,
\begin{align*}
\Bigl\lVert
  \alpha^\phi_a(Q_\epsilon)\phi_*(X)
& - \alpha^\phi_a(Q_\epsilon)
    \Bigl( \sum_{\{j: p(x_j) \le a\}} \phi(x_j)\phi(y_j)^* \Bigr)
  \Bigr\rVert \\
& \le \sum_{\{j: p(x_j) \not\le a\}} \norm{\alpha^\phi_a(Q_\epsilon)\phi(x_j)\phi(y_j)^*} \\
& < \epsilon \Bigl( \sum_{\{j: p(x_j) \not\le a\}} \norm{y_j} \Bigr).
\end{align*}
But using Lemma~\ref{lemma:alpharho} we have
\begin{align*}
\alpha^\phi_a(Q_\epsilon)
    \Bigl( \sum_{\{j: p(x_j) \le a\}} \phi(x_j)\phi(y_j)^* \Bigr)
& = \alpha^\phi_a(Q_\epsilon)\rho^\phi_a(1)
    \Bigl( \sum_{\{j: p(x_j) \le a\}}
    \rho^\phi_{p(x_j)}(\rankone{x_j}{y_j}) \Bigr) \\
& = \alpha^\phi_a(Q_\epsilon)
    \Bigl( \sum_{\{j: p(x_j) \le a\}}
    \rho^\phi_a((\rankone{x_j}{y_j}) \otimes 1) \Bigr) \\
& = \alpha^\phi_a(Q_\epsilon) \rho^\phi_a(T_a),
\end{align*}
and \cite[Proposition~1.12(2)]{fowrae} gives
$\norm{\alpha^\phi_a(Q_\epsilon) \rho^\phi_a(T_a)} = \norm{T_a}$,
so
\[
\norm{\alpha^\phi_a(Q_\epsilon)\phi_*(X)}
\ge \norm{T_a} - \epsilon \Bigl( \sum_{\{j: p(x_j) \not\le a\}} \norm{y_j} \Bigr).
\]
Since $\epsilon$ was arbitrary it follows that
\[
\norm{\phi_*(X)}
\ge \norm{\alpha^\phi_a(Q_\epsilon)\phi_*(X)} \ge \norm{T_a} = \norm X,
\]
and hence $\norm{\phi_*(X)} = \norm X$ as required.

(2) By Corollary~\ref{cor:wick}, finite sums of the form
$X := \sum_j i_E(x_j)i_E(y_j)^*$ are dense in $\cov(P,E)$,
so it suffices to fix such an $X$ and prove that
\begin{equation}\label{eq:contract}
\norm{\phi_*(\Phi_\delta(X))} \le \norm{\phi_*(X)}.
\end{equation}
It was shown in \cite[Proposition~5.5]{fowrae}
that there exists $a\in P$ such that
$\norm{\phi_*(\Phi_\delta(X))} = \norm{T_a}$,
where
\[
T_a := \sum_{\{j: p(x_j) = p(y_j) \le a\}} (\rankone{x_j}{y_j})\otimes 1^{p(x_j)^{-1}a}.
\]
We will show that for every $\epsilon>0$ there is
a nonzero projection $Q_\epsilon\in\Bb(\Hh)$ such that
\begin{equation}\label{eq:approx2}
\norm{\alpha^\phi_a(Q_\epsilon)\phi(x_j)\phi(y_j)^*\alpha^\phi_a(Q_\epsilon)} < \epsilon
\qquad
\text{unless $p(x_j) = p(y_j) \le a$.}
\end{equation}
Given this,
let $N$ be the number of $j$'s such that $p(x_j) = p(y_j) \le a$
does {\em not\/} hold, and proceed as in part (1) to estimate
\begin{align*}
\norm{\phi_*(X)}
& \ge \norm{\alpha^\phi_a(Q_\epsilon)\phi_*(X)\alpha^\phi_a(Q_\epsilon)} \\
& \ge \Bigl\lVert \alpha^\phi_a(Q_\epsilon)
          \Bigl( \sum_{\{j: p(x_j) = p(y_j) \le a\}} \phi(x_j)\phi(y_j)^* \Bigr)
          \Bigr\rVert - N\epsilon \\
& = \norm{\alpha^\phi_a(Q_\epsilon)\rho^\phi_a(T_a)} - N\epsilon \\
& = \norm{T_a} - N\epsilon \\
& = \norm{\phi_*(\Phi_\delta(X))} - N\epsilon.
\end{align*}
Letting $\epsilon$ decrease to zero gives \eqref{eq:contract}.

It therefore suffices to find a projection $Q_\epsilon$ satisfying \eqref{eq:approx2}.
For each $j$ such that $p(x_j)\le a$, $p(y_j) \le a$, $p(x_j) \ne p(y_j)$ and
$p(x_j)^{-1}a \vee p(y_j)^{-1}a < \infty$, define $d_j \in P\setminus \{e\}$
by
\[
d_j := \begin{cases}
(p(x_j)^{-1}a)^{-1}p(y_j)^{-1}a & \text{if $p(x_j)^{-1}a < p(y_j)^{-1}a$,} \\
(p(y_j)^{-1}a)^{-1}p(x_j)^{-1}a & \text{if $p(y_j)^{-1}a < p(x_j)^{-1}a$.} \end{cases}
\]
Fix $j$ such that $p(x_j) \le a$, $p(y_j) \le a$, and $p(x_j) \ne p(y_j)$.
If $Q\in\Bb(\Hh)$, then by Lemma~\ref{lemma:alpharho}(2) we have
\[
\alpha^\phi_a(Q)\phi(x_j)\phi(y_j)^*\alpha^\phi_a(Q)
= \phi(x_j)\alpha^\phi_{p(x_j)^{-1}a}(Q)\alpha^\phi_{p(y_j)^{-1}a}(Q)\phi(y_j)^*,
\]
which vanishes if $p(x_j)^{-1}a \vee p(y_j)^{-1}a = \infty$ since $\phi$ is covariant.
On the other hand, if $p(x_j)^{-1}a \vee p(y_j)^{-1}a < \infty$,
then $\alpha^\phi_{p(x_j)^{-1}a}(Q)\alpha^\phi_{p(y_j)^{-1}a}(Q)$ equals either
$\alpha^\phi_{p(x_j)^{-1}a}(Q\alpha^\phi_{d_j}(Q))$ or
$\alpha^\phi_{p(y_j)^{-1}a}(\alpha^\phi_{d_j}(Q)Q)$,
so $\alpha^\phi_a(Q)\phi(x_j)\phi(y_j)^*\alpha^\phi_a(Q)$
vanishes if $Q\alpha^\phi_{d_j}(Q) = 0$.
We conclude that any $Q_\epsilon$ which satisfies
\begin{equation}\label{eq:alphad}
Q_\epsilon\alpha^\phi_{d_j}(Q_\epsilon) = 0
\qquad\text{for each $j$,}
\end{equation}
and
\begin{equation}\label{eq:approx3}
\norm{\alpha^\phi_a(Q_\epsilon)\phi(x_j)\phi(y_j)^*\alpha^\phi_a(Q_\epsilon)} < \epsilon
\qquad
\text{if $p(x_j)\not\le a$ or $p(y_j)\not\le a$,}
\end{equation}
will also satisfy \eqref{eq:approx2}.

For the existence of such a projection $Q_\epsilon$, let
\[
F := \{x_j : p(x_j) \not\le a\} \cup \{y_j: p(y_j)\not\le a\},
\]
and let
\[
C := \{d_j\} \cup \{ a^{-1}p(x_j) : a < p(x_j)\} \cup \{ a^{-1}p(y_j) : a < p(y_j)\}.
\]
If the projection $R_C$ of Lemma~\ref{lemma:technical}(2) is nonzero,
we take it as $Q_\epsilon$;
the lemma assures that \eqref{eq:approx3} holds,
and since $R_C \le 1 - \alpha^\phi_{d_j}(1)$
and $\alpha^\phi_{d_j}(R_C) \le \alpha^\phi_{d_j}(1)$,
\eqref{eq:alphad} holds as well.
On the other hand, if $R_C = 0$, then for appropriately small $\delta > 0$
the unit vector $y_\delta$ provided by Lemma~\ref{lemma:technical}(2)
satisfies
\begin{equation}\label{eq:subprojection}
\norm{\alpha^\phi_a(\phi(y_\delta)\phi(y_\delta)^*)
      \phi(x_j)\phi(y_j)^*
      \alpha^\phi_a(\phi(y_\delta)\phi(y_\delta)^*)} < \epsilon
\end{equation}
if $p(x_j)\not\le a$ or $p(y_j)\not\le a$.
Lemma~\ref{lemma:aperiodic} then provides a unit vector $z$ such that
\[
\phi(y_\delta z)\phi(y_\delta z)^*
\alpha^\phi_{d_j}(\phi(y_\delta z)\phi(y_\delta z)^*) = 0
\qquad\text{for each $j$,}
\]
so $Q_\epsilon := \phi(y_\delta z)\phi(y_\delta z)^*$ satisfies \eqref{eq:alphad}.
Since $Q_\epsilon \le \phi(y_\delta)\phi(y_\delta)^*$,
we see from \eqref{eq:subprojection} that \eqref{eq:approx3}
holds as well, and the proof is complete.
\end{proof}

\begin{cor}\label{cor:simplicity2}
For each $\lambda$ belonging to some index set $\Lambda$,
let $(G^\lambda,P^\lambda)$ be a totally-ordered amenable group.
Let $E$ be a product system over $P:= *P^\lambda$.
If there exists $\mu\in\Lambda$ such that $E_s$ is infinite-dimensional
for every $s\in P^\mu\setminus\{e\}$, then $\cov(P,E)$ is simple.
\end{cor}

\begin{proof} By Proposition~\ref{prop:i_E covariant}, it suffices to show that
every covariant representation $\phi$ satisfies \eqref{eq:faithful}.
Suppose $s_1, \dots, s_n\in P\setminus\{e\}$ and $\dim E_{s_k} < \infty$ for every $k$.
Take any $s\in P^\mu\setminus\{e\}$, and fix $k$.
Then $s \not\le s_k$; for otherwise $E_{s_k} \cong E_s \otimes E_{s^{-1}s_k}$,
contradicting $\dim E_{s_k} < \infty = \dim E_s$.
On the other hand, $\dim E_{s_k}< \infty$ implies
$s_k \notin P^\mu\setminus\{e\}$, and hence $s_k \not\le s$.
From \eqref{eq:qto} we conclude that
$s\vee s_k = \infty$ for every $k$.  Since $\phi$ is covariant we thus have
$\rho^\phi_s(1)\prod_{k=1}^n (1 - \rho^\phi_{s_k}(1))
= \rho^\phi_s(1) \ne 0$,
and this implies \eqref{eq:faithful}.
\end{proof}

Taking $\abs\Lambda = 1$ gives the following generalization of \cite[Theorem~2.2]{dinhjfa}.

\begin{scholium}\label{scholium}
Let $(G,P)$ be a totally-ordered amenable group, and let $E$ be a product system over $P$.
If $E_s$ is infinite-dimensional for every $s\in P\setminus\{e\}$,
then $C^*(P,E)$ is simple.
\end{scholium}

\end{document}